\documentclass[10pt]{article}

\usepackage{amsmath}
\usepackage{amssymb}
\usepackage{amsthm}

\newcommand{\R}{\mathbb{R}}
\newcommand{\C}{\mathbb{C}}
\newcommand{\Z}{\mathbb{Z}}

\newcommand{\cc}[1]{\overline{#1}}
\renewcommand{\Re}{\mathop{\mathrm{Re}}}
\newcommand{\Tr}{\mathop{\mathrm{Tr}}}
\newcommand{\te}{\tilde}
\newcommand{\set}[2]{\left\{#1:#2\right\}}

\newtheorem{theorem}{Theorem}[section]
\newtheorem{lemma}[theorem]{Lemma}
\newtheorem{corollary}[theorem]{Corollary}

\newtheorem{proposition}[theorem]{Proposition}

\newtheorem{remark}{Remark}[section]

\numberwithin{equation}{section}

\begin{document}

\title{Rotation properties of isotropic dilation matrices}

\author{Victor G.~Zakharov\\[0.5ex]
\small Institute of Continuum Mechanics of Russian Academy of Sciences,\\
\small              Perm, 614013, Russia\\
\small  E-mail: \tt    victor@icmm.ru }

\maketitle{}

\begin{abstract}
In the paper, we consider bivariate isotropic dilation matrices
that are similar (up to constant factors) to rotation matrices;
and we show that, in this case,
the two-scale relation
can be considered also as
the relation between not only dilated but also rotated scaling functions.
We present sufficient conditions on a dilation matrix to be similar to a rotation matrix,
where the similarity transformation matrix is symmetric positive definite.
Also we present a simple test that
the dilation matrix performs the rotation by
an incommensurable to $\pi$ angle.
We show that if a dilation matrix is similar to a rotation matrix, then
an ellipse defined by the similarity transformation matrix
is invariant (up to a homogeneous dilation) under
the transformation by the dilation matrix.
\end{abstract}

\medskip

\noindent{\it Keywords:} \\[1ex]
Isotropic dilation matrices, similarity to rotation matrices,
rotation invariance, incommensurable to $\pi$ angles, scaling
functions, two-scale relation\\[1ex]


\section{Introduction}

Generally, the scaling function $\phi$ satisfies the {\em two-scale relation}
\begin{equation}\label{ScalingRelation}
  \phi(x) = \sum_{k\in{\Z\mathstrut}^2}h_{k}|\det A|^\frac{1}{2}\phi(Ax-k),
     \qquad x\in\R^2,
\end{equation}
where $A$ is called the {\em dilation matrix}.

Let us recall that a dilation matrix is called {\em isotropic} if the matrix
is diagonalizable over $\C$ and all its eigenvalues
are equal in absolute value. Any isotropic dilation matrix performs
an isotropic (homogeneous) expansion/contraction; and the isotropic dilation matrices have
some interesting properties,
see, for example,~\cite{R.-Q.Jia,Lindemann}.
Note that isotropic dilation matrices usually occur in the context of non-separable wavelet
bases.
(We can refer the reader to~\cite{Lindemann}, where
an extensive bibliography is presented.)

In the present paper, we show that, under
slightly more strong conditions, the isotropic real matrix $A$ is similar
(up to a constant factor)
to a rotation matrix
\begin{equation}\label{DilationMatrixSimilarity}
  A=(\det A)^{1/2}QRQ^{-1},
\end{equation}
where $R$ is the rotation matrix and the similarity
transformation matrix $Q$ is symmetric positive definite.
Thus such dilation matrix performs not only an expansion but also the rotation by some angle.
Consequently two-scale relation~\eqref{ScalingRelation}
can be also considered as a {\em two-angle} relation, i.\,e.,
as the relation between not only dilated but also rotated scaling functions.

In the paper, we have introduced a simple test
that the angle of the rotation performed by an isotropic dilation matrix
is incommensurable to $\pi$.
So if the corresponding scaling function/wavelet possesses some angular selectivity; then
the wavelet transform can detect {\em arbitrary} orientated features.

In general, the matrix $Q$ in~\eqref{DilationMatrixSimilarity}
distorts in some way the rotation performed by the matrix $R$. Nevertheless, in the paper,
we shall show that
an ellipse (defined by
$Q^{-2}$) will be invariant (up to a homogeneous dilation) under
the transformation by the dilation matrix $A$.
On the other hand, the invariance of the ellipse is equivalent to the invariance of the
corresponding quadratic form.
Note that we always can change the variables
(or, in other words, replace $\Z^2$ by another lattice in~\eqref{ScalingRelation})
such that the ellipse is transformed to a circle.
Thus, in this case, the dilation matrix will perform a pure rotation (and a homogeneous dilation).

Let us introduce here some general notation.

In the sequel,
we shall
not distinguish
vectors as points of the Euclidean space $\R^2$ and as column-matrices.
The Fourier transform of a function $f\in L^2(\R^2)\cap L^1(\R^2)$ is defined as
$$
  \hat f(\xi):=\int_{\R^2} f(x)e^{-i\xi\cdot x}\,dx,
$$
where $x:=(x_1,x_2)^T$, $\xi:=(\xi_1,\xi_2)^T$.
Two-scale relation~\eqref{ScalingRelation}
can be rewritten in the Fourier domain as
\begin{equation*}
  \hat\phi( \xi) = m_0\left(\left(A^*\right)^{-1} \xi\right)
     \hat\phi\left(\left(A^*\right)^{-1} \xi\right),     \qquad  \xi\in\R^2,
\end{equation*}
where the matrix $A^*$ is the Hermitian conjugate of the matrix $A$;
$m_0(\xi)$, $\xi\in\R^2$, is a $2\pi$-periodic function, which is called the {\em mask}.
By $R_\theta$ denote the (counterclockwise) rotation matrix by the angle $\theta$:
\begin{equation*}
  R_\theta:=\begin{pmatrix}
    \cos\theta & -\sin\theta \\
    \sin\theta & \cos\theta \\
  \end{pmatrix}.
\end{equation*}

\section{Similarity to rotation matrices}

\subsection{Theorem}

\begin{theorem}\label{Theorem_MyDecomposition}
Let $\te A$ be an $n\times n$ real {\em isotropic} matrix.
Let all the eigenvalues of the matrix $\te A$ be
equal to $1$ in absolute value,
then
\begin{equation}\label{MyDecomposition}
  \te A=QUQ^{-1},
\end{equation}
where $U$ is an orthogonal real
matrix and $Q$ is a symmetric positive definite
real matrix.
\end{theorem}

Here we present
the proof of Theorem~\ref{Theorem_MyDecomposition} for the bivariate case only.
To prove the theorem in the general case,
we can refer the reader to~\cite{Z_CA}.

\begin{proof}
Let $x:=\left(x_1,x_2\right)^T$ be an eigenvector of $\te A$ corresponding to an eigenvalue $\lambda$, i.\,e.,
$$
  \te A  x=\lambda  x;
$$
Complex conjugating both the sides of the previous expression,
we obtain
that $\cc{x}$ is the other
eigenvector of $\te A$ corresponding to the other eigenvalue $\cc\lambda$,
where the overline stands for the complex conjugation.
Since the matrix $\te A$ is isotropic; the matrix is diagonalizable:
\begin{equation}\label{MatrixADiagonalization}
  \te A =TDT^{-1},
\end{equation}
where $D=\left(%
   \begin{array}{cc}
     \lambda & 0 \\
     0 & \cc\lambda \\
   \end{array}%
  \right)$, $T=\left(%
    \begin{array}{cc}
      x_1 & \cc x_1 \\
      x_2 & \cc x_2 \\
    \end{array}%
  \right)$.
Using the polar decomposition, we can always present $T$ as follows
\begin{equation}\label{PolarDecomposition}
  T=QV,
\end{equation}
where $V$ is a unitary matrix and $Q$ is a positive definite Hermitian matrix:
$Q^2=TT^*$.
The matrix $Q^2$ can be written explicitly in component-wise form as follows
$$
  Q^2=\left(%
    \begin{array}{cc}
      x_1 & \cc x_1 \\
      x_2 & \cc x_2 \\
    \end{array}%
  \right)
  \left(%
    \begin{array}{cc}
      \cc x_1 & \cc x_2 \\
      x_1 & x_2 \\
    \end{array}%
  \right)
  =2\left(%
    \begin{array}{cc}
      |x_1|^2 & \Re(x_1\cc x_2) \\
      \Re(x_1\cc x_2) & |x_2|^2 \\
    \end{array}%
  \right).
$$
We see that $Q^2$ is a {\em real} matrix, thus
the ``square root'' $Q$ is also a real
matrix.
Using~\eqref{MatrixADiagonalization},~\eqref{PolarDecomposition}, the matrix $\te A$ can be written as follows
$$
  \te A = QVDV^{-1}Q^{-1}=QUQ^{-1},
$$
where $U:=VDV^{-1}=VDV^*$.
The matrix $U$ is a unitary matrix. Indeed,
$$
  UU^* = VDV^*VD^*V^*=V\left(%
     \begin{array}{cc}
       |\lambda|^2 & 0 \\
       0 & |\lambda|^2 \\
     \end{array}%
   \right)V^*=VIV^*=I,
$$
where $I$ is the identity matrix.

Finally, since the matrices $\te A$ and $Q$ in decomposition~\eqref{MyDecomposition}
are real, the matrix $U$ must be real (consequently orthogonal) also.
\end{proof}

\begin{remark}
The matrix $Q$ in formula~\eqref{MyDecomposition} is defined within a constant factor.
\end{remark}

\begin{remark}
Note that any $2\times2$ orthogonal matrix whose determinant is equal to 1 is
a rotation matrix. Thus,
in formulas like~\eqref{MyDecomposition},
if the determinant of the matrix $\te A$ is positive; then
we shall use the symbol $R$, which denotes a rotation matrix,
instead of $U$.
\end{remark}

\subsection{Two-angle relation}

\begin{proposition}
Let $A$ be a $2\times2$ real matrix with integer entries
and let $\det A>0$.
Let the matrix $\frac{1}{(\det A)^{1/2}}A$
satisfy the conditions of Theorem~\ref{Theorem_MyDecomposition}
and let $\theta\in[0,2\pi)$ be the angle of the rotation
performed by the orthogonal (actually rotation) matrix, which below we shall denote by $R$,
from~\eqref{MyDecomposition}.
Consider the set
\begin{equation*}
  \Theta:=\set{j\theta\bmod 2\pi}{j\in\mathbb{N}}.
\end{equation*}
Let $m_0$ be a mask and $\phi$ be the scaling function corresponding to $A$ and $m_0$;
then, for any angle $\vartheta\in\Theta$,
there exists an index
$j'\in\mathbb{N}$ such that
\begin{equation}\label{PhiInvarianceGeneralCase}
  \hat\phi(Q^{-1}R_{\vartheta}Q\, \xi)
    =\left(\prod_{j=1}^{j'}m_0((\det A)^{-j/2}Q^{-1}R^j Q\, \xi)\right)\hat\phi((\det A)^{-j'/2} \xi),
\end{equation}
\end{proposition}

If we change the frequency variables
\begin{equation}\label{CoordinateTransformXi}
  \xi\mapsto\xi':=Q\xi;
\end{equation}
then relation~\eqref{PhiInvarianceGeneralCase}
will have the simple form
\begin{equation}\label{PhiInvarianceAfterCoordinatTransform}
  \hat\phi(R_{\vartheta}\, \xi')
    =\left(\prod_{j=1}^{j'}m_0((\det A)^{-j/2}R^j\, \xi')\right)\hat\phi((\det A)^{-j'/2} \xi').
\end{equation}
The spatial variables transformation corresponding to~\eqref{CoordinateTransformXi}
is of the form
\begin{equation}\label{CoordinateTransformX}
  x\mapsto x':=Q^{-1}x.
\end{equation}
Using~\eqref{CoordinateTransformX},
formula~\eqref{PhiInvarianceGeneralCase} can be rewritten in the $x$-domain as follows
\begin{equation}\label{PhiInvarianceX-Domain}
  \phi(R_\vartheta  x') = \sum_{  k\in{\Z\mathstrut}^2}a_{  k}\phi((\det A)^{j'/2}  x'-  k),\qquad
          a_k\in\R.
\end{equation}

Thus formulas~\eqref{PhiInvarianceAfterCoordinatTransform},~\eqref{PhiInvarianceX-Domain}
(and~\eqref{PhiInvarianceGeneralCase})
can be interpreted as relations between not only scaled but also rotated scaling functions.

\begin{remark}
Note that the change of variables~\eqref{CoordinateTransformX}
(as well as~\eqref{CoordinateTransformXi}) can be interpreted as a coordinate transformation.
Moreover, map~\eqref{CoordinateTransformX}
can be considered also as the replacement of the lattice $\Z^2$
by another lattice $\Gamma:=Q\Z^2$ in two-scale relation~\eqref{ScalingRelation}.
\end{remark}

\begin{remark}
If the angle $\theta$ corresponding to the rotation matrix $R$
in formulas~\eqref{PhiInvarianceGeneralCase},~\eqref{PhiInvarianceAfterCoordinatTransform}
is incommensurable to $\pi$, then relation~\eqref{PhiInvarianceAfterCoordinatTransform}
(and~\eqref{PhiInvarianceGeneralCase},~\eqref{PhiInvarianceX-Domain})
is valid for a {\em dense} in $[0,2\pi)$ set of angles.
\end{remark}

\section{Rotation angle determination}

The following lemma allows determining the angle of rotation.
\begin{lemma}\label{AboutAngle}
Let a $2\times2$ real matrix $\te A$ be similar to a rotation matrix,
then the eigenvalues of $\te A$ are $\{e^{i\theta},e^{-i\theta}\}$,
where $\pm\theta$ is the rotation angle.
\end{lemma}
The proof is trivial.

\begin{remark}
By Lemma~\ref{AboutAngle}, the angle of rotation can be determined by the
eigenvalues up to an absolute value only, i.\,e., we cannot determine
the direction of the rotation.
\end{remark}

Let a dilation matrix $A$ be of the form
$$
  A=\begin{pmatrix}
    a & b \\
    c & d \\
  \end{pmatrix},\qquad a,b,c,d\in\Z.
$$
Suppose $\det A>0$ and
denote $\Delta:=\det A$, then $A$ can be presented as follows
$$
  A=\begin{pmatrix}
      \Delta^{1/2} & 0 \\
      0 & \Delta^{1/2} \\
    \end{pmatrix}\te A,
$$
where $\te A=\frac{1}{\Delta^{1/2}}
  \begin{pmatrix}
    a & b \\
    c & d \\
  \end{pmatrix}$.
The eigenvalues of $\te A$ can be explicitly calculated
$$
  \lambda_{1,2} = \frac1{2\Delta^{1/2}}\left(a+d\pm\sqrt{a^2-2ad+d^2+4bc}\right).
$$
Using the relation $\Delta=ad-bc$
and denoting $\Tr A:=a+d$ by $\tau$, we have
$$
  \lambda_{1,2} = \frac1{2\Delta^{1/2}}\left(\tau\pm\sqrt{\tau^2-4\Delta}\right).
$$

Consider three cases:
\begin{description}
\item[$\tau^2>4\Delta$:] $|\lambda_1|\ne|\lambda_2|$ ($\lambda_1,\lambda_2\in\R$);
\item[$\tau^2=4\Delta$:] $\lambda_1=\lambda_2$ ($\lambda_1,\lambda_2\in\R$);
\item[$\tau^2<4\Delta$:] $\lambda_{1,2}=\dfrac{\tau}{2\Delta^{1/2}}\pm
        i\dfrac{\sqrt{4\Delta-\tau^2}}{2\Delta^{1/2}}$.
\end{description}
The cases $\tau^2\ge4\Delta$ do not satisfy the conditions of
Theorem~\ref{Theorem_MyDecomposition}.
Consider the case $\tau^2<4\Delta$ in detail.

First we have $|\lambda_1|^2=|\lambda_2|^2=\dfrac{\tau^2}{4\Delta}
+\dfrac{4\Delta-\tau^2}{4\Delta}=1$. Thus the eigenvalues
can be presented as follows:
$\lambda_{1,2} =e^{\pm i\theta} =\cos\theta\pm i\sin\theta$,
where
\begin{equation}\label{CosSinTheta}
  \begin{aligned}
    \cos\theta=\dfrac{\tau}{2\Delta^{1/2}},\qquad
    \sin\theta=\dfrac{\sqrt{4\Delta-\tau^2}}{2\Delta^{1/2}}.
  \end{aligned}
\end{equation}
By Lemma~\ref{AboutAngle},
the corresponding rotation matrix is of the form
\begin{equation}\label{RotMatr}
  R :=
  \begin{pmatrix}
    \dfrac{\tau}{2\Delta^{1/2}} & \mp\dfrac{\sqrt{4\Delta-\tau^2}}{2\Delta^{1/2}}\\[1.7ex]
    \pm\dfrac{\sqrt{4\Delta-\tau^2}}{2\Delta^{1/2}} & \dfrac{\tau}{2\Delta^{1/2}}
  \end{pmatrix}.
\end{equation}

In Table~\ref{TableRTheta}, we present explicit forms of rotation matrices~\eqref{RotMatr} and
the corresponding rotation angles for some
$\tau$ ($\tau\ge0$) and $\Delta$.

\begin{table}[h!]
\begin{center}
\begin{tabular}{|c|c|c|}
\hline
$\tau,\Delta$ & $R$ & $\theta$ \\
\hline
 $\tau=0$: &
  $\begin{pmatrix}
    0 & \mp1 \\
    \pm1 & 0 \\
  \end{pmatrix}$ & $\pm\dfrac{\pi}{2}$ \\[2.ex]
\hline
$\tau^2=\Delta$: &
  $\begin{pmatrix}
    \frac12 & \mp\frac{\sqrt3}{2} \\
    \pm\frac{\sqrt3}{2} & \frac12 \\
   \end{pmatrix}$ & $\pm\dfrac{\pi}{3}$ \\[3.ex]
\hline
$\tau^2=2\Delta$: &
  $\begin{pmatrix}
    \frac{\sqrt2}{2} & \mp\frac{\sqrt2}{2} \\
    \pm\frac{\sqrt2}{2} & \frac{\sqrt2}{2} \\
   \end{pmatrix}$ & $\pm\dfrac{\pi}{4}$\\[3ex]
\hline
$\tau^2=3\Delta$: &
  $\begin{pmatrix}
    \frac{\sqrt3}{2} & \mp\frac12 \\
    \pm\frac12 & \frac{\sqrt3}{2} \\
   \end{pmatrix}$ & $\pm\dfrac{\pi}{6}$.\\[3ex]
\hline
\end{tabular}
\end{center}
\caption{Rotation matrices~\eqref{RotMatr} and the corresponding rotation angles $\theta$ for some
$\tau$ ($\tau\ge0$) and $\Delta$}\label{TableRTheta}
\end{table}

\begin{remark}
If the trace $\tau$ is negative, then
the rotation angle is of the form: $\theta+\pi$,
where $\theta$ is the angle
corresponding to the positive trace $|\tau|$.
\end{remark}

Note also that the case $\tau^2=4\Delta$ corresponds to the ``rotation'' matrix
$R=\begin{pmatrix}
    1 & 0 \\
    0 & 1 \\
  \end{pmatrix}$;
i.\,e., for this case, the rotation angle is zero.

Now we are interested in the rotation angles that do not correspond
to the cases presented in Table~\ref{TableRTheta}, i.\,e.,
\begin{equation}\label{TauConditions}
  \tau\in\Z:\qquad 0<\tau^2<4\Delta,\quad \tau^2\ne\Delta,2\Delta,3\Delta.
\end{equation}
It is rather surprise that,
for all cases~\eqref{TauConditions}, the angles of rotation are incommensurable to $\pi$.
This fact follows directly from the theorem of I.~Niven~\cite{Niven}.

\begin{theorem}[Ivan Niven, \cite{Niven}]\label{NivenTheorem}
If $\alpha/\pi$ and $\sin \alpha$ are both rational,
then the sine takes values $0$, $\pm1/2$, and $\pm1$.
\end{theorem}

Note that, since $\cos\alpha=\sin(\alpha-\pi/2)$;
I.~Niven's theorem is valid for the cosine also.

Consider $\cos2\theta$, where $\theta$ is defined by formulas~\eqref{CosSinTheta},
$$
  \cos2\theta = \cos^2\theta-\sin^2\theta = \dfrac{\tau^2-2\Delta}{2\Delta}.
$$
So $\cos2\theta$ is rational number. Consequently, by Theorem~\ref{NivenTheorem},
for all $\tau$ and $\Delta$ that satisfy~\eqref{TauConditions},
$\theta/\pi$ is an {\em irrational} number.

Summarize the results of this section in the following two statements.
\begin{lemma}\label{SummarizedTheorem1}
Let $A$ be a $2\times2$ real matrix with integer entries
and let $\det A>0$.
The matrix $\frac{1}{(\det A)^{1/2}}A$
satisfies the conditions of Theorem~\ref{Theorem_MyDecomposition} iff
$(\Tr A)^2<4\det A$.
\end{lemma}

\begin{theorem}\label{SummarizedTheorem2}
Let a matrix $A$ satisfy the condition of Lemma~\ref{SummarizedTheorem1}.
Suppose $(\Tr A)^2$ is equal to one of the following values:
$0$, $\det A$, $2\det A$, $3\det A$;
then the angle of the rotation performed by the corresponding rotation matrix
is $\pm\dfrac{\pi}{2}$, $\pm\dfrac{\pi}{3}$,
$\pm\dfrac{\pi}{4}$, $\pm\dfrac{\pi}{6}$, respectively; else
the rotation angle
is incommensurable to $\pi$.
\end{theorem}

In Table~\ref{Table},
we list the rotation angles for some values of the determinant and trace of
the matrix $A$.
(Recall that $(\Tr A)^2$ must not be more than $4\det A$.)
By Theorem~\ref{SummarizedTheorem2}, all the rotation angles listed in Table~\ref{Table}
that are not the explicit fractions $\dfrac{\pi}{\cdot}$
(but that are defined by the $\arccos(\cdot)$) are incommensurable to $\pi$.

\begin{table}[h]
\begin{tabular}{|c|c|c|c|c|c|}
\hline \rule[-0.2ex]{0pt}{1ex}
    & $\Tr A=0$ & $\Tr A=1$ & $\Tr A=2$ & $\Tr A=3$ & $\Tr A=4$ \\
\hline \rule[-2.2ex]{0pt}{5.5ex}
  $\det A=1$ & $\pm\dfrac{\pi}{2}$ & $\pm\dfrac{\pi}{3}$  &
0 & --- & --- \\
\hline \rule[-3.5ex]{0pt}{8.2ex}
  $\det A=2$ & $\pm\dfrac{\pi}{2}$ & $\pm\arccos\left(\dfrac{\sqrt2}{4}\right)$ &
$\pm\dfrac{\pi}{4}$ & --- & --- \\
\hline \rule[-3.5ex]{0pt}{8.2ex}
  $\det A=3$ & $\pm\dfrac{\pi}{2}$ & $\pm\arccos\left(\dfrac{\sqrt3}{6}\right)$ &
$\pm\arccos\left(\dfrac{\sqrt3}{3}\right)$ & $\pm\dfrac{\pi}{6}$ & --- \\
\hline \rule[-2.7ex]{0pt}{6.8ex}
  $\det A=4$ & $\pm\dfrac{\pi}{2}$ & $\pm\arccos\left(\dfrac{1}{4}\right)$ &
$\pm\dfrac{\pi}{3}$ & $\pm\arccos\left(\dfrac{3}{4}\right)$ & 0 \\
\hline \rule[-3.5ex]{0pt}{8.2ex}
  $\det A=5$ & $\pm\dfrac{\pi}{2}$ & $\pm\arccos\left(\dfrac{\sqrt5}{10}\right)$ &
$\pm\arccos\left(\dfrac{\sqrt5}{5}\right)$ & $\pm\arccos\left(\dfrac{3\sqrt5}{10}\right)$ &
   $\pm\arccos\left(\dfrac{2\sqrt5}{5}\right)$ \\
\hline
\end{tabular}
\caption{The rotation angles for some values of the determinant and trace of
a $2\times 2$ matrix $A$}\label{Table}
\end{table}

\begin{remark}
Note that, in the papers~\cite{GundyJonsson,LagariasWang},
the classification of the dilation matrices that have the
determinant values $\pm2$ and are similar to unimodular matrices
has been presented. It is interesting that the line of
Table~\ref{Table} for $\det A=2$ corresponds to the matrix
classification in~\cite{GundyJonsson,LagariasWang}. However any
connection with the rotation properties of the dilation matrices
has not been discussed in the
papers~\cite{GundyJonsson,LagariasWang}.
\end{remark}

\section{Invariance of ellipse shape}

\begin{corollary}\label{CorollaryFromMyDecomposition}
Under the conditions of Theorem~\ref{Theorem_MyDecomposition}, we have
\begin{align}
    \label{Q-2Invariance}
    &\te A^TQ^{-2}\te A=\left(\te A^T\right)^{-1}Q^{-2}\te A^{-1}=Q^{-2},\\
    \label{Q2Invariance}
    &\te AQ^2\te A^T=\te A^{-1}Q^2\left(\te A^T\right)^{-1}=Q^2.
\end{align}
\end{corollary}
Using~\eqref{MyDecomposition}, the proof is straightforward.

Consider a quadratic form
\begin{equation}\label{QuadraticForm}
  W(x):={x}^TQ^{-2}x,\quad x\in\R^2,
\end{equation}
where $Q$ is the similarity transformation matrix from~\eqref{MyDecomposition}.
Recall that the matrix $Q$ (and, consequently, $Q^{-2}$) is positive definite; thus
$W(x)$ is a positive definite quadratic form.
By Corollary~\ref{CorollaryFromMyDecomposition}, we see that quadratic
form~\eqref{QuadraticForm} is invariant (up to a constant factor)
under the map $x\mapsto x':=Ax$.
Indeed, using~\eqref{Q-2Invariance}, we have
\begin{equation*}
  W(x')=W(Ax)=x^TA^{T} Q^{-2}Ax
      = \det A\,x^TQ^{-2}x = \det A\, W(x).
\end{equation*}
Thus the transformation by
the matrix $A$ will save the shape of the ellipse
\begin{equation}\label{Ellipse}
  {x}^TQ^{-2}x=C,\qquad x\in\R^2,
\end{equation}
where $C>0$ is a constant.
Note that, using~\eqref{CoordinateTransformX},
quadratic form~\eqref{QuadraticForm} gets a paraboloid of revolution $W(x)=|x|^2$,
$x\in\R^2$; and
ellipse~\eqref{Ellipse} turns into a circle.

Similarly, the quadratic form
${x}^TQ^{2}x$, $x\in\R^2$,
is invariant under the transformation
by the matrix $A^T$, see~\eqref{Q2Invariance}.
Note that the axes of the ellipse ${x}^TQ^{2}x=C$ are
orthogonal to the corresponding axes of ellipse~\eqref{Ellipse}.

\begin{remark}
We suppose that, multiplying by
an appropriate real value, the matrix $Q^{-2}$ (as well as $Q^{2}$)
that corresponds to the matrix $A$ (see~\eqref{MyDecomposition})
with integer entries
can be made a matrix with integer entries also. This will be discussed elsewhere.
\end{remark}

\section{Examples}

As a very simple example we consider a (quincunx) dilation matrix
$\begin{pmatrix}
    1 & -1 \\
    1 & 1 \\
  \end{pmatrix}$.
Obviously, this matrix performs the rotation by angle $\pi/4$
and the isotropic dilation by $\sqrt2$, cf.\ Table~\ref{Table};
and the corresponding matrix $Q$ is the identity matrix.
Thus the invariant ellipse is the circle. (Note that the well-known quincunx
dilation matrix is symmetric and differs from the matrix above.)

Consider another example
\begin{equation}\label{MatrixExample}
  A
  =\begin{pmatrix}
    0 & -2 \\
    1 & 1 \\
  \end{pmatrix}.
\end{equation}
We see that $\det A=2$ and $\Tr A=1$.
By Table~\ref{Table}, it follows that the matrix performs rotation by the angle that is
incommensurable to $\pi$: $\arccos\left(\frac{\sqrt2}{4}\right)\approx
69{.}2951889^\circ$.
The corresponding matrix $Q^{-2}$ is of the form
$$
  Q^{-2}=
  \begin{pmatrix}
    1 & \frac12 \\
    \frac12 & 2 \\
  \end{pmatrix}
$$
and quadratic form~\eqref{QuadraticForm} is $W(x,y):=x^2+xy+2y^2$.
Note that matrix~\eqref{MatrixExample} has been also considered in the
papers~\cite{GundyJonsson,LagariasWang}.
However rotation properties, in particular,
that the rotation angle is incommensurable to $\pi$ have not been noted.

\section{Conclusion}

Using the ellipse (quadratic form) invariance, we can choose a mask
and construct a (compactly supported) scaling function that represents polynomials
from the null-space of the elliptic operator corresponding to the invariant quadratic form.
We refer the reader to~\cite{Z_CA} for details.
Moreover, these scaling functions (and the corresponding wavelets) possess some interesting
angular properties; and the scaling functions can be considered as {\em compactly supported} counterparts
of the polyharmonic B-splines~\cite{Rabut}, see also~\cite{MicchelliRabutUtreras,VilleBluUnser}.
This is the object of the forthcoming paper~\cite{Z_IJWMIP}.

\section*{Acknowledgements}
We would like to thank the founder of the Wikipedia Jimmy~D.~Wales and all other contributors
who help us to find a reference to the statement
to complete our paper.

Research was partially supported by RFBR grant No.~12-01-00608-a.

\end{document}